\documentclass[a4paper,12pt]{article}

\usepackage{latexsym}
\usepackage{amssymb}
\usepackage{amsthm}

\theoremstyle{definition}

\renewcommand{\L}{\ensuremath{\mathrm{L}}}
\newcommand{\Con}{\ensuremath{\mathrm{C}}}
\newcommand{\Cinf}{\ensuremath{\mathrm{C}^\infty}}
\newcommand{\D}{\ensuremath{{\cal D}}}

\newcommand{\E}{\ensuremath{{\cal E}}}

\newcommand{\mb}[1]{\ensuremath{\mathbb{#1}}}
\newcommand{\N}{\mb{N}}

\newcommand{\R}{\mb{R}}

\newcommand{\cl}[1]{\ensuremath{\mathrm{cl}[#1]}}

\newcommand{\G}{\ensuremath{{\cal G}}}
\newcommand{\Gsg}{\ensuremath{{\cal G}_\mathrm{s,g}}}

\newcommand{\EM}{\ensuremath{{\cal E}_{\mathrm{M}}}}

\newcommand{\NN}{\ensuremath{{\cal N}}}

\renewcommand{\d}{\ensuremath{\partial}}
\newcommand{\rd}{\ensuremath{\tilde{\partial}}} 

\newcommand{\curl}{\ensuremath{\mbox{\rm curl}\,}}

\renewcommand{\div}{\ensuremath{\mbox{\rm div}\,}}

\newcommand{\au}{\frac{u}{\sqrt{1+u^2}}}
\newcommand{\wu}{\sqrt{1+u^2}}

\newfont{\bl}{msbm10 scaled \magstep2}

\newcommand{\bem}[1]{\vadjust{\rlap{\kern\hsize\thinspace\vbox%
                       to0pt{\hbox{${}_\clubsuit${\small\tt #1}}\vss}}}}

\newcommand{\map}{\ensuremath{\rightarrow}}

\newcommand{\col}{\colon}

\newcommand{\notmid}{\mid\kern-0.5em\not\kern0.5em}

\newcommand{\norm}[2]{{\| #1 \|}_{#2}}
\newcommand{\lone}[1]{\norm{#1}{\L^1}}

\newcommand{\linf}[1]{\norm{#1}{\L^\infty}}

\newcommand{\al}{\alpha}

\newcommand{\de}{\delta}
\newcommand{\eps}{\varepsilon}

\newcommand{\vphi}{\varphi}

\newcommand{\sig}{\sigma}

\newcommand{\supp}{\mathop{\mathrm{supp}}}

\renewcommand{\em}[1]{{\it #1\/}}
\newcommand{\ovl}[1]{\overline{#1}}

\begin{document}
\begin{center}
{\Large \bf Regularized derivatives in a 2-dimensional model of 
  self-interacting fields with singular data}

\vspace{3mm}

{\bf G\"unther H\"ormann and Michael Kunzinger}
\end{center}

\vspace{1cm}

\paragraph{Abstract.}
The coupled Maxwell-Lorentz system describes feed-back action of electromagnetic 
fields in classical electrodynamics. When applied to point-charge sources (viewed as
limiting cases of charged fluids) the resulting nonlinear weakly hyperbolic system 
lies beyond the scope of classical distribution theory. Using regularized derivatives
in the framework of Colombeau algebras of generalized functions we analyze a 
two-dimensional analogue of the Maxwell-Lorentz system. After establishing existence and 
uniqueness of solutions in this setting we derive some results on distributional limits
of solutions with $\delta$-like initial values.

\vspace{1mm}
{\bf Keywords:} {\it Regularized derivatives, Colombeau algebras, 
                     self-interaction, weakly hyperbolic, Maxwell-Lorentz}

{\bf AMS subject classification:} Primary 35Dxx, 46F10 , secondary 35Lxx, 78A25

\vspace{4mm}

{\large \bf 1. Introduction}\setcounter{section}{1}\setcounter{equation}{0}

\vspace{2mm}

In recent years there has been considerable development in the application of
nonlinear theories of generalized functions to problems in physical field 
theory, most successfully in general relativity (cf.\ \cite{cvw}, \cite{stein},
\cite{stein2}, \cite{wilson})
but also starting in electrodynamics (see \cite{novi96}). 
Due to the nonlinearities inherent in the coupling of fields and 
equations of motion methods from classical distribution 
theory soon turn out to be insufficient for treating such problems.  
On the other hand, the theory of algebras of generalized functions in the
sense of J.F.~Colombeau offers the possibility of applying large classes
of nonlinear operations to distributional objects. 

We shortly recall the principal setup of electrodynamic field theory.
Configuration space is modeled by 4-dimensional Minkowski space. The basic 
physical 
quantities are represented by the charge- and current density $J_0$, $J$, 
the electric field $E$, and the magnetic field $B$. Maxwell's field equations are
\[
\begin{array}{clccl}
\mbox{(H1)} & \div B = 0 & & \mbox{(H2)} & \d_t B + \curl E = 0\\
\mbox{(I1)} & \div E = J_0& & \mbox{(I2)} & \d_t E - \curl B = -J\\
\end{array}
\]
with initial conditions
\[
 E|_{t=0} = E_0\qquad  B|_{t=0} = B_0\; .
\]
\hrulefill \hspace{8cm}\ \\
Institut f\"ur Mathematik, Universit\"at Wien,  A-1090 Wien, Austria\\
Email: guenther.hoermann@univie.ac.at, michael.kunzinger@univie.ac.at

These equations imply the following well-known compatibility
relations 
\[
\begin{array}{lcc}
\mbox{(C0)} & & \div B_0 = 0,\; \div E_0 = J_0|_{t=0}\\
\mbox{(C1)} & & \d_t J_0 + \div J = 0\; .
\end{array}
\]
On the other hand (C0), (C1) together with (H2), (I2), and the
initial conditions are sufficient to derive equations
(H1), (I1) (take the divergence of equations (H2), (I2) and integrate
with respect to time). So, in order to solve the complete Maxwell system
it is sufficient to study the hyperbolic Cauchy problem consisting of initial
data and equations (H2), (I2) while $E_0$, $B_0$, $J$, and
$J_0$ are assumed to obey the relations (C0), (C1).

Let us now investigate the nature and role of the charge density $J_0$
and the current $J$. In an idealized model we will interpret these as
quantities of a charged fluid (or sometimes clouds of charged dust, 
cf.\ \cite{parrott}) with charge density proportional to mass density
(for convenience we set this ratio equal to $1$).
In this model we write 
\begin{equation}\label{fluid}
(J_0,J) = \rho (u_0,u)\; ,
\end{equation}
where $\rho$ represents charge density distribution in space-time and
$(u_0,u)$ is the 4-velocity field of the fluid, normalized with respect to
the Minkowski metric by
\begin{equation}\label{vel_norm}
u_0^2 - |u|^2 = 1\; , 
\end{equation}
which is plausible as the 4-velocity is always assumed to be time-like.

If we study the feed-back action of the electromagnetic field on the
charged fluid the equations of motion for $u$ are determined
by the Lorentz force according to (cf.\cite{parrott})
\[
  \d_t u + \sum_{j=1}^3 \frac{u_j}{\sqrt{1+|u|^2}} \d_j u = 
  E + \frac{u}{\sqrt{1+|u|^2}}\times B \; . 
\]

Observe that the above formulae involve products of the basic physical
quantities (resp.\ of their coefficients). These quantities can
be singular if for example jump discontinuities of certain
conductor-dielectrica configurations are to be modeled or if point charges
are studied.

For instance, if one is interested in the equations of motion of charged
particles resulting from the electromagnetic field they generate
then $\rho$ -- and consequently $E$ and $B$ -- involve strong singularities
corresponding to concentration of charge on the particles' world line.

Therefore in considering the realistic problem of self-interaction we have
to deal with nonlinear differential equations involving singular functions and
products and compositions thereof.

Even worse, the equations of motion cannot be written down consistently 
in the framework of classical distribution theory. Therefore, often physical 
motivations are given for renormalization and adding corrective terms.
Unfortunately these ``tricks'' give rise to strange behaviour of 
solutions to the resulting equations (e.g.\ the Lorentz-Dirac
equation). Critical discussions of the emerging inconsistencies 
and the resulting unphysical solutions (``run-away-solution'', 
``preaccelaration'', \ldots) can be found in the books \cite{thirring} and
\cite{parrott}.

The aim of this paper is to analyze the mathematical structure of 
the coupled Maxwell-Lorentz system by means of a low dimensional ``toy model'' 
which displays some of the main features of the general problem. The plan
of exposition is as follows: In section 2 we develop a (1+1)-dimensional model
of self-interacting charges which will provide the basis of our further investigations.
As the model equations are nonlinear, of evolution type and involve singular initial
data we use the formalism of regularized derivatives in Colombeau algebras of generalized 
functions (cf. \cite{coheimo}). Section 3 presents the main existence and uniqueness
result in this framework.  Finally, in section 4 we consider singular initial data 
(point charges). We illustrate the influence of the regularization procedure on the 
existence of distributional limits of the unique solutions established above.


\vspace{4mm}

{\large \bf 2. A (1+1)-dimensional model of self-interacting
charges}\setcounter{section}{2}\setcounter{equation}{0}

\vspace{2mm}

We try to mimic the purely mathematical structure of Maxwell equations 
with only one space dimension, i.e.\ in a two-dimensional Minkowski space
$M$ with Lorentz metric in coordinates $(t,x)$ equal to $dt^2-dx^2$. 
The field quantities (electric and magnetic) are scalar functions on $M$ and 
will be denoted by $E$ and $B$.  
Further ``physical'' quantities are represented by the scalar charge density 
$\rho$ and the relativistic (normalized) velocity field 
$(\sqrt{1+u^2},u)$ of a ``fluid'' consisting of a collection of 
(strictly time-like world lines of) charged particles. 
Therefore $\rho (\sqrt{1+u^2},u)$ represents the relativistic current of the 
fluid. Interpreting divergence as simple $x$-derivative and removing curl
terms in the original Maxwell equations this produces the electromagnetic fields 
according to
\begin{eqnarray*}
\d_x B = 0 && \d_t B = 0\\
\d_t E = -\rho u && \d_x E = \rho \sqrt{1+u^2} \; .
\end{eqnarray*}
Therefore we may set $B$ equal to some constant $B_0$.
The second line defines a gradient condition on $E$ which is
solvable if and only if the following (relativistic) charge conservation 
holds
\begin{equation}
\d_t \big(\rho\wu\big) + \d_x \big(\rho u\big) = 0 \; .\label{Ladung}
\end{equation}
In order to recover an evolution type system we simply combine the
equations for $E$ into one by adding them.

Finally the feed-back action of the field on the charged particles of the fluid 
is modeled by mimicking the structure of the Lorentz force law:
\begin{equation}
\d_t u + \au\d_x u = E + \au B \, .\label{Lorentz}   
\end{equation}
The left hand side evaluated along a world line $z$ with
$\dot{z}(s)=u(z(s))$ gives $\ddot{z}(s)$, i.e.\ an actual acceleration.

Again, the above formulae involve products (and even more complicated
nonlinear operations) of the basic physical
quantities. 
By the transformation $\rho \mapsto \sig = \rho \wu$ we can rewrite the 
system of equations in evolution form 
\begin{equation}
 \d_t\pmatrix{E\cr u\cr \sig} +
 \pmatrix{ 1 &  0   & 0\cr
           0 &  \au & 0\cr
           0 &  \frac{\sig}{(1+u^2)^{3/2}} & \au
         } \d_x\pmatrix{E\cr u\cr \sig} =
 \pmatrix{  \sig-\sig\au \cr E+\au B_0 \cr 0} \; . \label{system}  
\end{equation}
Recall that $E$, $u$, and $\sig$ are real valued functions of
the variables $(t,x)$ and $B_0$ is a real constant. We study the
Cauchy problem for system (\ref{system}) with initial data
\begin{equation}\label{initial}
  \pmatrix{E\cr u\cr \sig}  =  \pmatrix{E_0\cr u_0\cr \sig_0} \ \ \
  \mbox{ when } t=0. 
\end{equation}

System (\ref{system}) is quasi-linear hyperbolic but not strictly hyperbolic. Also,
its coefficients are not uniformly bounded.  Furthermore, the main application 
we have in mind is to model self-interaction of charged particles 
(approximated by fluid clusters shrinking to isolated points) which 
qualitatively corresponds to $\sig_0$ being proportional to a delta distribution.
In order to address such questions we need a mathematical theory which is capable of
describing singular objects and allows for an unrestricted application of
both differentiation and of a sufficiently large class of nonlinear operations. 
Differential algebras of generalized functions containing the space of distributions
as a subspace and the space of smooth functions as a faithful subalgebra and providing
the above-mentioned tools were introduced by J.F. Colombeau in \cite{col1} (see also \cite{col2},
\cite{mo}). We will therefore work in this framework.

Using the relation $\d_x(1+u^2)^{1/2} = u(1+u^2)^{-1/2}\d_x u$ and denoting
$a(y) = y(1+y^2)^{-1/2}$ we can rewrite the above system in the more  
compact form 
\begin{eqnarray}
\d_t E + \d_x E &=& \sig\, (1-a(u)) \label{E_eqn}\\
\d_t u + \d_x (1+u^2)^{\frac{1}{2}} &=& E + B_0\, a(u) \label{u_eqn}\\
\d_t \sig + \d_x (\sig a(u)) &=& 0 \label{sig_eqn} \; . 
\end{eqnarray}                

This form of the equations will be our starting point for the method
of regularized derivatives which is the subject of the following section.


\vspace{4mm}

{\large \bf 3. Regularized derivatives -- existence and uniqueness of
   solutions in Colombeau algebras} \setcounter{section}{3}\setcounter{equation}{0}

\vspace{2mm}

We will work with a variant of Colombeau algebras
denoted by $\Gsg(\ovl{\Omega})$ ($\Omega$ an open subset of $\R^n$) which 
is a simplified version with representatives globally bounded and smooth
(cf. \cite{coheimo}):
denote by $\D_{\L^\infty}(\ovl{\Omega})$ the algebra of restrictions of
smooth functions on $\R^n$ to $\ovl{\Omega}$ which have all derivatives
in $\L^\infty(\R^n)$; then within the algebra $\E_\mathrm{s,g}(\ovl{\Omega}) = 
\D_{\L^\infty}(\ovl{\Omega})^{(0,\infty)}$ we define
\begin{eqnarray*}
\E_\mathrm{M,s,g}(\ovl{\Omega}) =  
    \{ (u_\eps)_\eps\in\E_\mathrm{s,g}(\ovl{\Omega}) \mid 
    \forall\al\in\N_0^n\, \exists p>0 : 
    \norm{\d^\al u_\eps}{\L^\infty(\ovl{\Omega})}\! = O(\eps^{-p})\,(\eps\to 0) \} \\
\NN_\mathrm{s,g}(\ovl{\Omega}) =      
    \{ (u_\eps)_\eps\in\E_\mathrm{s,g}(\ovl{\Omega}) \mid 
    \forall\al\in\N_0^n\, \forall q>0 : 
    \norm{\d^\al u_\eps}{\L^\infty(\ovl{\Omega})}\! = O(\eps^{q})\,(\eps\to 0) \}
    \, ;
\end{eqnarray*}
then $\Gsg(\ovl{\Omega})$ is the quotient algebra of 
$\E_\mathrm{M,s,g}(\ovl{\Omega})$ modulo $\NN_\mathrm{s,g}(\ovl{\Omega})$.
In our application we will have $\Omega = (-T,T)\times\R$ for $T>0$ arbitrary.

The notion of regularized derivative was introduced in \cite{coheimoI}-\cite{coheimo}. 
It extends weak solution concepts in
algebras of generalized functions. The idea is to not only represent
singular functions by (classes of) regularizations but also to consider
differentiation as a limiting case of more smooth operations, namely convolution
with derivatives of delta regularizations. This means that a partial
derivative $\d_j U$ of the generalized function $U=\cl{(u_\eps)_\eps}$
will be replaced by $\cl{(\d_j \psi_\eps * u_\eps)_\eps}$ where
$(\psi_\eps)_\eps$ is a delta net (in view of application to our model only 
first order derivatives are considered here). Thus if $U$ is associated to
a distribution $T$ the regularized derivative will be associated to $\d_j T$.

Regularized derivatives provide a powerful solution concept for a wide variety
of partial differential equations. For example, the Cauchy problem for quasilinear evolution
type systems 
\begin{eqnarray*}
&&\d_t u(t,x) = \sum_{|\alpha|\le m} A_\alpha(u(t,x))\d^\alpha u(t,x) + B(u(t,x))\\
&& u(0,x) = u_0(x)
\end{eqnarray*}
is uniquely solvable in this framework provided that $A_\alpha$ and $B$ are bounded in all
derivatives (\cite{coheimo}, Thm. 5.1).

The choice of which of the occurring derivatives are to be regularized may 
depend on special information about the model, the structure and form of the
equations, and the type of effects one wants to observe. In our case we will 
systematically replace all $x$-derivatives in the model equations of the form
(\ref{E_eqn}-\ref{sig_eqn}) but leave the time derivatives unchanged in order 
to reflect the evolution type of the equations. Note that this is different 
from regularizing $x$-derivatives in (\ref{system}).

The regularized $x$-derivative is defined by choosing a function 
$\vphi\in\D$ with integral equal to $1$ on $\R$ --
we call this a \em{mollifier} henceforth -- and an 
increasing 
function $h\col (0,\infty)\map (0,\infty)$ with $h(\eps)\to 0$ as $\eps\to 0$. 
Setting $\vphi_\nu(x) = \vphi(x/\nu)/\nu$ we have $\d_x (\vphi_\nu(x)) = 
\vphi_\nu'(x) = \vphi'(x/\nu)/\nu^2$ (by $\vphi_\nu'$ we will always understand
$(\vphi_\nu)'$) and define the regularized derivative of a
generalized function $U = \cl{(u_\eps)_\eps}$ by the formula
\begin{equation}\label{reg_der}
  (\rd_x)_h U = \cl{(\vphi_{h(\eps)}' * u_\eps)_\eps} \; .
\end{equation}
The function $h$ can be used to control the speed of 
convergence towards the usual derivative compared to the regularization
parameter $\eps$. That it is well-defined and also the basic properties
of the regularized derivative follow from the estimates
\begin{eqnarray}
  \norm{\d^\al\d_x \vphi_{h(\eps)}*u_\eps}{\L^\infty} \leq
    \norm{\vphi}{\L^1} \norm{\d^\al\d_x u_\eps}{\L^\infty} \\
  \norm{\d_x \vphi_{h(\eps)}*u_\eps}{\L^\infty} \leq
    h(\eps)^{-1} \norm{\vphi'}{\L^1} \norm{u_\eps}{\L^\infty} \, .
\end{eqnarray}

Now we can state the main existence and uniqueness result for the 
Cauchy problem with regularized $x$-derivatives. We will only have to impose 
the following growth condition on $h$:

$\forall p\in \N \ \exists k\in \N \ \mbox{ s.t.} \ \exp(h(\eps)^{-p}) = 
O(\eps^{-k}) \ \ \ (\eps\to 0)\hspace*{\fill} (*)$

For example, $h(\eps) = C \mbox{ln}(|\mbox{ln}(\eps)|)^{-1}$ satisfies
$(*)$.
\paragraph{Theorem 1.} Let $T>0$ and assume that $h$ satisfies $(*)$.
Then given $E_0$, $u_0$, $\sig_0$ in $\Gsg(\R)$ the system
\begin{eqnarray}
\d_t E + (\rd_x)_h E &=& \sig\, (1-a(u)) \label{E_rd}\\
\d_t u + (\rd_x)_h (1+u^2)^{\frac{1}{2}} &=& E + B_0\, a(u) \label{u_rd}\\
\d_t \sig + (\rd_x)_h (\sig a(u)) &=& 0 \label{sig_rd}  \; ,
\end{eqnarray}  
with initial conditions 
\begin{equation}\label{initial_G}
  E\mid_{_{t=0}} = E_0 ,\quad u\mid_{_{t=0}} = u_0 ,\quad 
  \sig\mid_{_{t=0}} = \sig_0  
\end{equation}
has a unique solution $(E,u,\sig)$ in $\Gsg([-T,T]\times\R)^3$.
\begin{proof} We will proceed in two steps: 
First, the candidate for a Colombeau solution will be built up by nets of
smooth solutions at fixed $\eps$. Second, existence and uniqueness in \Gsg
will be proved by verifying additional estimates with respect to the 
regularization-parameter $\eps$. To begin with, we consider global
solvability of the corresponding integro-differential equations where we
will suppress all references to the index $\eps$.

\em{Step 1:} integrating with respect to $t$ and inserting initial conditions
we get an equivalent system of integral equations
\begin{eqnarray}
  \pmatrix{E \cr u \cr \sig}\!\!(t,x) &=& \pmatrix{E_0 \cr u_0 \cr \sig_0}\!\!(x)-
  \int\limits_0^t \pmatrix{E(\tau,.)*\vphi' \cr
    \big(1+u(\tau,.)^2\big)^{1/2}*\vphi' \cr
     a(u(\tau,.))\sig(\tau,.)*\vphi'}\!\!(x)\, d\tau + \nonumber\\
 && + \int\limits_0^t \pmatrix{\sig(\tau,x)\big(1-a(u(\tau,x))\big) \cr
    E(\tau,x) + B_0 a(u(\tau,x)) \cr 0} d\tau \, . \label{int_equ}
\end{eqnarray}
Writing $V=(v_1,v_2,v_3)$ instead of $(E,u,\sig)$ and denoting 
$$ \norm{V}{T} = \max_{j=1,2,3} \sup_{|t|\leq T,x\in\R}|v_j(t,x)|
$$
and 
$V_0=(E_0,u_0,\sig_0)$, which is globally bounded by assumption, we
consider the right hand side as an operator $R$ on the set
\begin{equation}\label{BT}
  B_T(V_0) = \{ V\in\Con([-T,T]\times\R)^3 \mid \norm{V-V_0}{T}\leq 1 \}\, .
\end{equation}
Using $\norm{a}{\L^\infty}\leq 1$, $\sqrt{1+s^2}\leq 1+s$ for $s\geq 0$,
and Young's inequality for the convolutions, one derives from (\ref{int_equ}) 
\[
 |(RV)_j(t,x)-(V_0)_j(x)| \leq T \norm{\vphi'}{\L^1}(1+\norm{V}{T}) +
   T (2\norm{V}{T}+|B_0|) \quad (j=1,2,3)
\]
which shows that $RV\in B_T(V_0)$ for $T$ small enough. To see that $R$
is actually a contraction on $B_T(V_0)$ for small $T$ we can estimate
the terms occurring in $\norm{RV-RW}{T}$ for $V$, $W\in B_T(V_0)$ after 
applying Young's inequality
and/or taking $\sup$ for all integrands. First, in the integral 
involving convolutions the non-trivial factors are handled as follows:  
by the mean value theorem
applied to the function $f(y)=\sqrt{1+y^2}$ and $\norm{f'}{\L^\infty} =
\norm{a}{\L^\infty}\leq 1$ we have $\sup|\sqrt{1+v_2^2}-\sqrt{1+w_2^2}| 
\leq \norm{V-W}{T}$; since also $a'$ is globally bounded by $1$ and by 
definition of $B_T(V_0)$ we further 
estimate 
\begin{eqnarray*}
  |a(v_2)v_3-a(w_2)w_3| \leq |a(v_2)(v_3-w_3)|+|\big(a(v_2)-a(w_2)\big)w_3| \\
  \leq |v_3-w_3|+|v_2-w_3||w_3| \leq \norm{V-W}{T}(1+\norm{V_0}{\L^\infty})
  \, ;
\end{eqnarray*}
For dominating the $\sup$ of the integrands in the second integral we simply 
observe
\begin{eqnarray*}
|v_3\big(1-a(v_2)\big) - w_3\big(1-a(w_2)\big)| \leq |v_3-w_3| + 
   |a(v_2)v_3-a(w_2)w_3| \\
|v_1+B_0a(v_2) - w_1-B_0a(w_2)| \leq |v_1-w_1| + |B_0||a(v_2)-a(w_2)|
\end{eqnarray*}
to get terms we already estimated above. Altogether for
$T$ small, depending on $|B_0|$, $\norm{V_0}{\L^\infty}$, and 
$\norm{\vphi'}{\L^1}$ the map $R\col B_T(V_0)\to B_T(V_0)$ is a contraction. 
Therefore a solution of the fixed point equation exists for small $T$.

Assuming that $V$ and $W$ are two fixed points for $T>0$ arbitrary one can 
estimate their difference similar to the above except that only $\sup_{x\in\R}$ 
is to be considered inside integrals with respect to $\tau$. A standard Gronwall 
argument then yields global uniqueness of the solution.  An a priori
estimate for a solution $V$ in $[-T,T]\times\R$ is derived by using the fixed 
point formula (\ref{int_equ}), $|(1+s^2)^{1/2}-1| \leq |s|$, taking $\sup_{x}$
inside the integrals, and again applying Gronwall's lemma:
\[
  \norm{V}{T} \leq \big(\norm{V_0}{\L^\infty}+
     T(\norm{\vphi'}{\L^1}+|B_0|)\big) e^{T(\norm{\vphi'}{\L^1}+1)} \; .
\] 
Therefore existence of a global solution follows. Clearly this solution
is bounded on $[-T,T]\times\R$ by definition of $B_T(V_0)$. 

Finally to get unique solutions in $\D_{\L^\infty}([-T,T]\times\R)$ from 
initial data in $\D_{\L^\infty}(\R)$ we proceed inductively. Set 
$B_T^0(V_0)=B_T(V_0)$ with $\norm{V}{T,0}=\norm{V}{T}$; then for $k\in\N$ we 
define
$$
  B_T^k(V_0)=\{V\in\Con^k([-T,T]\times\R)^3 \mid \norm{V-V_0}{T,k}\leq 1\}
$$
where 
$$
  \norm{V}{T,k} = \norm{V}{T,k-1} + 
    T \big(\norm{\d_t V}{T,k-1} + \norm{\d_x V}{T,k-1}\big) \; .
$$
Then estimates can be done essentially by the same arguments as above 
because the nonlinear combinations of unknowns remain of the same
type (e.g.\ all derivatives of the function $a$ are again bounded
and therefore application of Young's inequality, mean value theorem,
and Gronwall's lemma are still possible). The weight factor $T^k$
for derivatives of order $k$ in the above norms serves to control 
inner derivatives and assures again e.g.\ the contraction property
for small $T$. This concludes the purely classical base for the next step.

\em{Step 2:} Existence and uniqueness of solutions in
$\Gsg([-T,T]\times\R)^3$.\\
In equation (\ref{int_equ}) all involved functions except $\vphi$ 
now carry an index $\eps$ and $\vphi$ is to be replaced by $\vphi_{h(\eps)}$.
We assume $((E_{0\eps},u_{0\eps},\sig_{0\eps}))_\eps\in\E_\mathrm{M,s,g}(\R)^3$.
In order to establish existence of a solution we will show that the net 
$(V_\eps)=((E_\eps,u_\eps,\sig_\eps))_\eps$ belongs to 
$\E_\mathrm{M,s,g}([-T,T]\times\R)^3$.

From the integral equations for $(V_\eps)$ and using Gronwall-type arguments as
above we conclude
\begin{eqnarray*}
&&\max_{1\le i\le 3}\sup_{|t|\leq T,x\in\R}|v_{i,\eps}(t,x)| \le \\
&& \le C(\linf{V_{0,\eps}},T,\lone{\vphi'})|B_{0,\eps}|h(\eps)^{-1}
\exp(3Th(\eps)^{-1}\lone{\vphi'}),
\end{eqnarray*}                         
which (by $(*)$) gives the \EM-estimates of order $0$.     
The \EM-estimates for the $x$-derivatives of higher order follow inductively by
similar arguments if in each step we write the $x$-derivative of any
convolution product of $\vphi^{(m)}$ with any function $f$ as $\vphi^{(m+1)}*f$
and use $(*)$. The \EM-estimates for the $t$-derivatives follow directly by
induction from the differential equations.

Concerning uniqueness, suppose that $W = \cl{(W_\eps)_\eps}$ is another solution in
$\Gsg([-T,T]\times\R)^3$. This means that there exists some
$N=\cl{(N_\eps)_\eps} \in\NN([-T,T]$ $\times$ $\R)^3$ with
\begin{eqnarray*}
\d_t  \pmatrix{v_{1,\eps}-w_{1,\eps} \cr v_{2,\eps}-w_{2,\eps} \cr v_{3,\eps}-w_{3,\eps}}
   &=& -\pmatrix{(v_{1,\eps}-w_{1,\eps}) \cr 
   [(1-v_{2,\eps}^2)^{\frac{1}{2}}-(1-w_{2,\eps}^2)^{\frac{1}{2}}]\cr
    [a(v_{2,\eps})v_{3,\eps}-a(w_{2,\eps})w_{3,\eps}]}*\vphi'_{h(\eps)} + \\
   && \pmatrix{v_{3,\eps}(1-a(v_{2,\eps})) - w_{3,\eps}(1-a(v_{3,\eps}))\cr 
   v_{1,\eps}-w_{1,\eps}+B_{0,\eps}(a(v_{2,\eps})-a(w_{2,\eps})) \cr 0} 
   + \pmatrix{N_{1,\eps}\cr N_{2,\eps}\cr N_{3,\eps}}
\end{eqnarray*}
From this by arguments similar to those used above we inductively derive
\NN-estimates for each derivative of $(V_\eps-W_\eps)_\eps$. Thus $V=W$ in
$\Gsg([-T,T]\times\R)^3$.

\end{proof}

\paragraph{Remark 1.}
Having obtained a unique solution to the field equations let us reinvestigate
the behaviour of point particles subjected to the velocity field $(\sqrt{1+u^2},u)$,
i.e. determine its integral curves in $\R^2$.
Classically, this amounts to considering the system
\begin{eqnarray}
&& \dot z_0 (s) = \sqrt{1 + u(z(s))^2} \nonumber \\
&& \dot z_1 (s) = u(z(s)) \label{ode}\\
&& z(0) = (t_0,x_0) \nonumber
\end{eqnarray} 
where $z = (z_0,z_1):I\to\R^2$ ($I$ some interval containing $0$).
By the first line of (\ref{ode}) we may reparametrize via $r = z_0(s)$, thereby obtaining
the following equation for $(r,w(r)) = z(s)$:
\begin{eqnarray} 
&& \dot w(r) = a(u(r,w(r))) \label{ode2}\\
&& w(t_0) = x_0 \nonumber
\end{eqnarray} 
Of particular interest
to us is $u$ arising as a solution of system (\ref{E_rd}) - (\ref{sig_rd}) with strongly
singular (e.g. $\delta$-like) initial data. In this case, (\ref{ode2}) is an example of 
a nonlinear ODE involving generalized functions. Note that its right-hand side contains a composition
of generalized functions (defined by componentwise insertion of representatives).  
A theory of such equations in 
the framework of Colombeau generalized functions has been developed by R. Hermann and M. Oberguggenberger
in \cite{hermo,hermo2}. Since all derivatives of $a$ are globally bounded and $u \in \Gsg$ it follows that
$\cl{(a\circ u_\eps)_\eps}$ is in $\Gsg$ and therefore of $L^\infty$-type (see \cite{hermo2}, Def. 2.3).
Thus \cite{hermo2}, Thm. 3.3 ensures the existence of a solution in $\G$ which is actually
in $\Gsg$ (which follows inductively from (\ref{ode2})). In particular, this allows to derive
(generalized) trajectories of charged particles without having to resort to renormalization procedures 
altering the actual equations of motion. This may indicate a way of avoiding the inconsistencies
mentioned at the end of section 1.


\vspace{4mm}

{\large \bf 4. Point charges as initial data}\setcounter{section}{4} \setcounter{equation}{0}

\vspace{2mm}

As we mentioned in the introduction we will now study the situation
where the initial data $(E_0,u_0,\sig_0)$ are qualitatively of the
form $(0,0,\de)$. This can be modeled by choosing $\sig_0$ to be the
class of a \em{strict delta net} $(\rho^\eps)_\eps$ in $\Gsg(\R)$ in the
sense of \cite{mo}, Def.7.1. According to this definition $(\rho^\eps)_{\eps>0}$ 
is a family of test functions in $\D(\R)$ such that 
\begin{eqnarray}\label{strict_delta}
&\supp \rho^\eps \to \{ 0 \} \quad\mathrm{as\ } \eps\to 0\nonumber &\\
&\int\! \rho^\eps(x) \, dx = 1 \quad\mathrm{for\ all\ } \eps >0 &\\
&\int\! |\rho^\eps(x)| \, dx \quad\mathrm{is\ bounded\ independently\ of\ } 
  \eps\;. &  \nonumber
\end{eqnarray} 
By the above theorem we obtain a unique solution in $\Gsg([-T,T]\times\R)^3$ 
of system (\ref{E_rd}-\ref{sig_rd}) with (\ref{initial_G}) for appropriate
scaling function $h$. However, the dependence of the definition of the
regularized derivative on the mollifier $\vphi$ is essential.
This was observed already in \cite{cohei} and illustrated by numerical
tests. Also in the context of scalar conservation laws it was shown how
the choice of ``left'' or ``right'' mollifiers allows to recover the classical
entropic solutions. 

We are not able to give a complete description of conditions for the existence
of distributional limits of solutions of the more complicated system (\ref{E_rd}-\ref{sig_rd}).
Nevertheless we will present some results analyzing the consequences of employing  
one-sided mollifiers.

\paragraph{Proposition 1.} Let $T>0$, $h$ as in Thm.1, 
and $(E_0,u_0,\sig_0) = (0,0,\cl{(\rho^\eps)_\eps})$ with $(\rho^\eps)_\eps$
a strict delta net. If $\supp \vphi \subseteq [-\infty,0]$ (resp. $\subseteq [0,\infty]$) then  
the unique solution $(E,u,\sig)\in\Gsg([-T,T]\times\R)$ to 
(\ref{E_rd}-\ref{sig_rd}) with initial conditions (\ref{initial_G}) is
supported in $[-T,T]\times(-\infty,0]$ (resp. $[-T,T] \times [0,\infty)$).
\begin{proof} It will suffice to consider the case of a left-sided mollifier (since a
similar argument applies to the case of a right-sided mollifier). 
Fix some $x_0>0$. By our assumption on the support of $\vphi$ we
have
\[
\sig_\eps(t,x) = - \frac{1}{h(\eps)} \int\limits_0^t \int\limits_{-1}^0
\vphi'(y) a(u_\eps(\tau,x-h(\eps)y))\sig_\eps(\tau,x-h(\eps)y) dy d\tau
\]
for $x\ge x_0$ and $\eps$ sufficiently small (depending on $x_0$). 
Also, since $a$ is globally bounded 
\[
\sup_{x\ge x_0} |\sig_\eps(t,x)| \le \frac{C}{h(\eps)} \lone{\vphi'} 
\int\limits_0^t \sup_{x\ge x_0} |\sig_\eps(t,x)| d\tau,
\]
so $\sup_{x\ge x_0} |\sig_\eps(\tau,x)| = 0$ by Gronwall's inequality which
(since $x_0$  was arbitrary) establishes the claim for $\sigma$. 
Then since 
\[
\sup_{x\ge x_0} |E_\eps(t,x)| \le \sup_{x\ge x_0}|\sig_\eps(t,x)(1-a(u_\eps(t,x)))| + 
\frac{C}{h(\eps)} \lone{\vphi'} 
\int\limits_0^t \sup_{x\ge x_0} |E_\eps(t,x)| d\tau,
\]
it follows that
$\sup_{x\ge x_0} |E_\eps(t,x)| = 0$. Finally, 
\begin{eqnarray*}
&& \d_x u_\eps(t,x) = 
- \frac{1}{h(\eps)} \int\limits_0^t \int\limits_{-1}^0
\vphi'(y) a(u_\eps(\tau,x-h(\eps)y))\d_x u_\eps(\tau,x-h(\eps)y) dy d\tau + \\
&& \hphantom{\d_x u_\eps(t,x) =}
+ B_0 \int\limits_0^t a'(u_\eps(\tau,x)) \d_x u_\eps(\tau,x) d\tau
\end{eqnarray*}
From this it follows that $\d_x u_\eps(t,x)=0$ for $x>0$. But then since
$(1+u_\eps^2)^{1/2}*\vphi_{h(\eps)}' = (a(u_\eps)\d_x u_\eps)*\vphi_{h(\eps)} =
0$ in this region we obtain
\[
\d_t u_\eps(t,x) = B_0 a(u_\eps(t,x)) = B_0 \int\limits_0^t a'(u_\eps(\tau,x))
\d_t u_\eps(\tau,x) d\tau \,.
\]
Thus also $\d_t u_\eps(t,x)=0$ and therefore $u_\eps(t,x)=0$ for $x>0$.
\end{proof}

Applying this result we now demonstrate how the choice of a 
left mollifier can destroy the intuitive expectation of field propagation 
along the characteristics of the original system.  

\paragraph{Corollary 1.} The (unique) solution of (\ref{E_rd}-\ref{sig_rd}), (\ref{initial_G})
 according to Prop.1 cannot be associated to a distribution in 
 $\D'(\R^2)^3$.
\begin{proof} 
Assume that $(E,u,\sig)\in \D'(\R)^3$ is the distributional limit 
of $(E_\eps,u_\eps,\sig_\eps)$. Then from $\d_t E_\eps +
\vphi_{h(\eps)}'*E_\eps = \sig_\eps(1-a(u_\eps))$ it follows that $\sig_\eps
a(u_\eps) \to \sig - F$, where $F = \d_t E + \d_x E$. Also, the third line of
the regularized system implies $\d_t\sig + \d_x (\sig -F) = 0$, so
\begin{equation}
\d_t (\sig -\d_x E) + \d_x (\sig -\d_x E) = 0 \label{sigE}
\end{equation}
Using the coordinate transformation $T(t,x) = (x,t-x)$ it follows that $\sig - \d_x E$
is of the form $v(t-x)$ (i.e. $T^*(1\otimes v)$) for some $v\in \D'(\R)$, hence
is contained in $\Cinf(\R,\D'(\R))$. 
By Prop.1 the supports of $\sig$ and $E$ are contained in $x\le
0$ and from the initial conditions we get $(\sig - \d_x E)(0)= \delta(x)$. Thus
$\sig - \d_x E = \delta(t-x)$, which contradicts $\supp(\sig - \d_x E)\subseteq
\{x\le 0\}$. 
\end{proof}

It remains an open question whether there are choices of the mollifier 
$\vphi$ which yield distributional limits of the solution for singular
initial data as in Prop.1.


\paragraph{Remark 2.}
Let us interpret equations (\ref{E_eqn}-\ref{sig_eqn}) in the form
$P(E,u,\sig)=0$ with the operator $P:\Cinf(\R^2)^3\to\Cinf(\R^2)^3$
given by
\begin{equation}\label{P}
  P(E,u,\sig) = \pmatrix{
        \d_t E + \d_x E + \sig\big(a(u)-1\big)\cr
        \d_t u + a(u) \d_x u - E - B_0 a(u)\cr
        \d_t\sig - \d_x\big(a(u)\sig\big) 
  }\; .
\end{equation}
We can decompose $P$ into a sum of a linear differential operator $L$ and 
a (purely) nonlinear part $Q$ with
\begin{eqnarray}
L(E,u,\sig) = \pmatrix{
        \d_t E + \d_x E - \sig \cr
        \d_t u - E \cr
        \d_t\sig  
  } \label{L}\\
Q(E,u,\sig) = \pmatrix{
        a(u)\sig \cr
        a(u) \big(\d_x u - B_0\big) \cr
        - \d_x\big(a(u)\sig\big) 
  }\; . \label{Q}
\end{eqnarray} 

We observe that in a linearization of $P$ (in the sense of \cite{rauch}, p.20) 
at a point $(t_0,x_0)$ with $u(t_0,x_0)=0$ all terms corresponding to $Q$
vanish and we are left exactly with $L$. Therefore, if we consider some smooth
approximation of the initial data (with $q$ the constant point charge)
\begin{equation}\label{init}
  (E_0,u_0,\sig_0) = (0,0,q \de)
\end{equation}
and a linearization of the system near $t=0$, the unique distributional solution
$(E,u,\sigma)\in C^1([0,T],\D'(\R))^3$ of the Cauchy problem
\begin{equation}\label{L_sys}
    \pmatrix{
        \d_t E + \d_x E \cr
        \d_t u \cr
        \d_t\sig  
  }         =   \pmatrix{\sig\cr E \cr 0}
\end{equation}
with initial condition (\ref{init}) may give a hint on
how to construct approximate smooth solutions of the non-linear system. 
The solution to the linearized system is given by 
\begin{equation}\label{lin_sol}
\pmatrix{E\cr u\cr \sig}(t,x) = q \cdot \pmatrix{H(x) - H(x-t) \cr
                                       (t-x)\big(H(x) - H(x-t)\big) \cr
                                       \de(x)} \; ,
\end{equation}
where $H$ denotes the Heaviside function. However, using a standard regularization
of this solution (by means of a model delta net) as a candidate for such an approximate
solution runs into difficulties: explicit calculations show that the blow up of the 
term $\sig a(u)$ at $x=0$ is unavoidable unless regularizations with $\delta H \approx 0$
are considered.

\end{document}